\magnification 1200

\hyphenation{di-men-sio-nal}

\font\Bbb=msbm10
\def\BBB#1{\hbox{\Bbb#1}}

\def\R{\BBB R}
\def\C{\BBB C}
\def\Z{\BBB Z}
\def\N{\BBB N}

\def\Ci{{\cal C}^\infty}
\def\K{{\cal K}}
\def\RR{{\cal F}}

\def\Diff{\hbox{\it Diff}}
\def\Map{\hbox{\it Map}}
\def\Tr{\hbox{\it Tr}}
\def\diver{\hbox{\it div}}

\def\TN{{\BBB T}^N}
\def\GLN{GL_N}
\def\GL{GL}
\def\glN{gl_N}
\def\gl{gl}
\def\slN{sl_N}

\def\al{\alpha}
\def\be{\beta}

\def\d{\partial}
\def\x{{\bf x}}
\def\v{{\bf v}}
\def\w{{\bf w}}

\def\tf{{\widetilde f}}
\def\tg{{\widetilde g}}

\centerline{\bf Abelian Extensions of the Group of Diffeomorphisms of a 
Torus.}

\

\centerline {\bf Yuly Billig}

\centerline{School of Mathematics and Statistics, Carleton University,}
\centerline{1125 Colonel By Drive, Ottawa, K1S 5B6, Canada}
\centerline{e-mail: billig@math.carleton.ca} 

\

{\bf Abstract.}
In this Letter we construct abelian extensions of the group of diffeomorphisms
of a torus. We consider the jacobian map, which is a crossed homomorphism from
the group of diffeomorphisms into a toroidal gauge group. A pull-back under this 
map of a central 2-cocycle on a gauge group turns out to be an abelian cocycle
on the group of diffeomorphisms. We show that in the case of a circle, the 
Virasoro-Bott cocycle is a pull-back of the Heisenberg cocycle. We also give
an abelian generalization of the Virasoro-Bott cocycle to the case of a manifold
with a volume form.

\


\


\

{\bf 0. Introduction.}

\

The Lie group that corresponds to the Virasoro Lie algebra is the central
extension of the group of diffeomorphisms of a circle with the Virasoro-Bott
cocycle (2.1). It is well-known that the Virasoro cocycle on the Lie
algebra of vector fields on a circle does not admit a generalization
as a central 2-cocycle to the case of higher dimensional manifolds.
Nonetheless there exist abelian extensions of the Lie algebra of vector fields
on a torus that generalize the Virasoro Lie algebra. These abelian extensions 
play an important role in the representation theory of toroidal Lie algebras
([EM], [L], [BB], [Bi]).

The goal of this Letter is to construct the corresponding abelian extensions for
the group of diffeomorphisms of a torus.

We approach this problem by linking the abelian cocycles on the group of 
diffeomorphisms of a torus with the central cocycles on the toroidal gauge groups.
The central extensions of the toroidal gauge groups were studied in [LMNS] (see
also [PS], [FK]).

The connection between the group of diffeomorphisms of a torus and toroidal gauge
groups is given by the jacobian map. Let us briefly outline this correspondence.
The differential of a diffeomorphism of a torus is a mapping of the tangent bundle
into itself. However, using the fact that torus has a trivial tangent bundle, we can 
globally identify all tangent spaces with $\R^N$. After this identification, the
differential of a diffeomorphism of a torus becomes a mapping of a torus into 
$\GLN(\R)$. In the standard coordinates, this mapping is given by the jacobian
matrix. In this way we obtain the jacobian map from the group of diffeomorphisms
of a torus into a toroidal gauge group:
$$ J: \Diff(\TN) \rightarrow \Map(\TN, \GLN (\R)) .$$
In all constructions of this Letter the torus may be replaced with an arbitrary 
manifold with a trivial tangent bundle.

The jacobian mapping is not a group homomorphism, but what is called a crossed
homomorphism (see definition in Section 3). In case $N=1$, we get a crossed 
homomorphism from the group of diffeomorphisms of a circle into abelian loop
group $\Map(S^1, \R^*)$. A central extension of this loop group is an infinite 
dimensional Heisenberg group. We point out that the Virasoro-Bott cocycle may be
interpreted as a pull-back of the Heisenberg cocycle under the jacobian map
(Theorem 2.1). This observation motivates our constructions of 2-cocycles on the
group of diffeomorphisms in higher dimensions.

First we describe our construction of 2-cocycles by pull-back in a very
general set-up. We show (Lemma 3.1) that under a crossed homomorphism of two
groups 
$$ j: D \rightarrow M,$$
a pull-back of a $D$ invariant central 2-cocycle on $M$ is an abelian 2-cocycle
on $D$.

We apply this lemma to construct an abelian extension of the group of
diffeomorphisms of a manifold with a volume form (Theorem 4.1), generalizing 
the Virasoro-Bott group, for which the manifold is the circle. This extension
(as well as the Virasoro-Bott extension) trivializes on the subgroup of volume
preserving diffeomorphisms.

We also construct another extension of the group of diffeomorphisms of a torus 
(Theorem 3.2) that remains non-trivial when restricted to the subgroup of volume
preserving diffeomorphisms. 

We would like to mention here that Ovsienko and Roger described abelian extensions 
of the group of diffeomorphisms of a circle with the modules of tensor densities [OR]. 

The structure of the Letter is the following. In Section 1 we present a geometric
approach to the abelian extensions of the Lie algebra of vector fields on a torus
and introduce the jacobian map from the group of diffeomorphisms of a torus into 
a toroidal gauge group. In Section 2 we give two constructions of the 
infinite dimensional Heisenberg group and exhibit the link between the Virasoro-Bott
cocycle and the Heisenberg cocycle. In Section 3 we describe the general 
construction of a pull-back of a cocycle under a crossed homomorphism and use this
approach to get an abelian extension of the group of diffeomorphisms of a torus.
Finally in Section 4 we obtain another, inequivalent abelian cocycle that generalizes
the Virasoro-Bott cocycle to the case of a manifold with a volume form.

\

\

{\bf 1. Abelian extensions of the Lie algebra of vector fields of a torus.}

\

It is well-known that the Lie algebra $W_N$ of vector fields on an 
$N$-{}dimensional torus $\TN$ (Witt algebra) has a non-trivial central 
extension only when $N=1$.
In that case there exists a unique non-trivial central extension of $W_1$, 
which 
is the Virasoro Lie algebra. The Virasoro cocycle can not be generalized as 
a central 2-cocycle to higher dimensions. However when $N > 1$, there exist 
{\it abelian} extensions of $W_N$ that generalize the Virasoro cocycle. 
These abelian extensions appeared in the representation theory [EM], [L],
where it was discovered that it is easier to construct representations for 
certain abelian extensions of $W_N$ rather than for the Witt algebra itself.

Let us describe the construction of these abelian extensions of $W_N$. 

Let $\RR$ be the algebra of $\Ci$ real-valued functions on a torus 
$\TN = \R^N / \Z^N$. Note that the algebra of complex-valued functions on a 
torus, $\C \otimes \RR$ has a convenient topological basis
$$ \left\{ \exp \left( 2\pi i \sum\limits_{j=1}^N r_j x_j \right) 
\, \bigg| \, r_j 
\in 
\Z \right\} .$$
Often this basis is used in algebraic setting (considering Fourier 
polynomials instead of Fourier series). However, the group of 
diffeomorphisms of a torus can not be realized using Fourier polynomials, 
so we will be working with $\Ci$ functions instead. We define the group
$\Diff(\TN)$ as the group of $\Ci$ diffeomorphisms of a torus.
An element $F \in \Diff(\TN)$ can be viewed as an $N$-tuple of functions
$F = (F_1(\x), \ldots, F_N(\x) ),$ such that the map
$$ \x = (x_1, \ldots, x_N) \mapsto (F_1(\x), \ldots, F_N(\x) )$$
is a bijection $\TN \rightarrow \TN$.

The Lie algebra $W_N$ consists of the vector fields $\v = 
\sum\limits_{j=1}^N v_j (\x) {\d \over \d x_j} $ with
$v_j (\x) \in \RR$.

Consider the space of differential $k$-forms on $\TN$:
$$ \Omega^k = \left\{ \sum\limits_{1 \leq j_1 < \ldots < j_k \leq N}
a_{j_1\ldots j_k} (\x) \, dx_{j_1} \wedge \ldots \wedge dx_{j_k} 
\; \bigg| \;
a_{j_1\ldots j_k} (\x) \in \RR \right\} .$$
There is a differential map $d: \Omega^k \rightarrow \Omega^{k+1}$:
$$ d \left( a(\x) dx_{j_1} \wedge \ldots \wedge dx_{j_k} \right) = 
\sum_{j=1}^k {\d a \over \d x_j} dx_j \wedge dx_{j_1} \wedge \ldots \wedge 
dx_{j_k} .$$

The Lie algebra of vector fields acts in an obvious way on $\RR$:
$ \v \cdot a(\x) = \sum\limits_{j=1}^N v_j (\x) {\d a \over \d x_j}$. 

The algebra of functions $\RR$ is also a right module for the group $\Diff(\TN)$
$$ a(\x) (F) = a(F(\x)), \quad a \in \RR, \, F\in\Diff(\TN) ,$$
and so is the spaces of differential forms:
$$ \left( a(\x) dx_{j_1} \wedge \ldots \wedge dx_{j_k} \right) (F) =
a(F(\x)) dF_{j_1}(\x) \wedge \ldots \wedge dF_{j_k} (\x). \eqno{(1.1)}$$
It is easy to see that the differential map is a homomorphism of $\Diff(\TN)$ 
modules.

The action of the group $\Diff(\TN)$ on $k$-forms gives rise to the action of its
Lie algebra $W_N$:
$$ \v \cdot \left( a(\x) dx_{j_1} \wedge \ldots \wedge dx_{j_k} \right)  =$$
$$(\v \cdot a(\x)) dx_{j_1} \wedge \ldots \wedge dx_{j_k} + 
\sum\limits_{r=1}^k
a(\x) dx_{j_1} (\x) \wedge \ldots \wedge dv_{j_r}(\x) \wedge \ldots
\wedge dx_{j_k} .$$
This action is called the Lie derivative action.

Next we need to discuss the jacobian of a diffeomorphism of a torus, and the 
jacobian of a vector field.

Let $F \in \Diff(\TN)$, $F: \TN \rightarrow \TN$. The differential of $F$ is a 
map 
of tangent bundles:
$$ dF: T(\TN) \rightarrow T(\TN), $$
with linear maps 
$$T_{\x} (\TN) \rightarrow T_{F(\x)} (\TN) .$$
However torus is a flat manifold and its tangent bundle is trivial:
$T(\TN) = \TN \times \R^N$, so all tangent spaces can be globally 
identified with $\R^N$, and we define the jacobian of $F$ as a map
$$ F^J: \TN \rightarrow Hom(\R^N, \R^N) .$$
Since $F$ is an invertible transformation of the torus, $F^J$ is an 
invertible linear homomorphism, so 
$$ F^J: \TN \rightarrow \GLN (\R) .$$
In coordinates $(x_1, \ldots x_N)$, the matrix of the jacobian of a 
diffeomorphism $F$ is 
$$F^J = \left( {\d F_i \over \d x_k} \right)_{i,j = 1, \ldots, N}.$$

The jacobian map $J: F \mapsto F^J$ associates with every element of $\Diff(\TN)$
an element of a gauge group $\Map (\TN, \GLN)$. However the resulting map 
$$J: \Diff(\TN) \rightarrow \Map (\TN, \GLN) \eqno{(1.2)}$$ 
is not a group homomorphism, 
but what is called a {\it crossed homomorphism} [K]:
$$ (FG)^J = F^J (G) \cdot G^J \quad \hbox{\rm (chain rule)}. \eqno{(1.3)}$$
In the above formula we used the action of
$\Diff(\TN)$ on the toroidal gauge group 
\break
$\Map (\TN, \GLN)$
by automorphisms:
$$ m(\x) \mapsto m(F(\x)), \hbox{\rm \ where \ }
m \in \Map (\TN, \GLN), \quad F \in \Diff(\TN), \quad \x \in \TN . 
\eqno{(1.4)}$$

In a similar way, the jacobian of a vector field $\v$ is defined as a 
matrix
$$\v^J = \left( {\d v_i \over \d x_k} \right)_{i,j = 1,
\ldots, N},$$
thus giving a linear map from $W_N$ to the toroidal analog $\Map (\TN, \glN) 
\cong \glN \otimes \RR$ of the loop algebra.

Now we can use the jacobians to give a simple descriptions of the abelian 
extensions of $W_N$. We will actually consider two possible extension spaces.
In the first case, the 2-cocycles will have values in the factor module 
$\K = \Omega^1 / d \Omega^0$ . 
This module is of special interest since it is the space of the universal 
central extension for the toroidal Lie algebras.
Two such $\K$ valued 
2-cocycles have appeared in the representation theory ([EM], [L], [BB]):
$$ \tau_1 (\v, \w) = \Tr \left( \v^J d \w^J \right) $$
and 
$$ \tau_2 (\v, \w) = \Tr \left( \v^J \right)  \Tr \left( d\w^J \right). $$
Viewed as $\Omega^1$ valued expressions these are not 2-cocycles, but modulo
$d \Omega^0$ these become skew-symmetric and satisfy the Jacobi identity.

There is a closely related pair of 2-cocycles on $W_N$, obtained from $\tau_1$,
$\tau_2$ by applying the differential map $d: \Omega^1 \rightarrow \Omega^2$.
Since $d^2 = 0$, this map factors through the factor module $\Omega^1 / d \Omega^0$.
The differentials $d \tau_1, d \tau_2$ are thus $\Omega^2$ valued 2-cocycles
on $W_N$:
$$ d \tau_1 (\v, \w) = \Tr \left( d \v^J \wedge d \w^J \right) $$
and
$$ d \tau_2 (\v, \w) =  \Tr \left( d\v^J \right) \wedge  \Tr \left( d\w^J 
\right). $$
Dzhumadildaev proved that for $N>1$ the space of $\Omega^2$ valued 2-cocycles 
$H^2(W_N, \Omega^2)$ is two-{}dimensional and is spanned by $d\tau_1, d\tau_2$ 
[D].

These cocycles give 2-parametric abelian extensions of the Lie algebra of vector
fields on a torus with extension spaces $\K$ or $\Omega^2$ (see e.g. Section 
2 of [BB]).

The main goal of this Letter is to construct the corresponding  
abelian extensions for the group of diffeomorphisms of a torus.

\

{\bf Remark 1.} We can see that linear combinations of $\tau_1, \tau_2$ are in 1-1
correspondence with the symmetric invariant bilinear forms on $\glN (\R)$. Since
for $N>1$ the Lie algebra $\glN (\R)$ is not simple, but only 
reductive, $\glN (\R) = 
\slN (\R) \oplus \R I$, the space of symmetric invariant forms on $\glN (\R)$
is two-{}dimensional and is spanned by the forms $\Tr (AB)$ and $\Tr(A) \Tr(B)$.
Thus a linear combination of $\tau_1$ and $\tau_2$ can be written as
$$\tau (\v, \w) = (\v^J | d\w^J),$$
where $(\cdot | \cdot)$ is a symmetric invariant bilinear form of $\glN 
(\R)$.  

\

{\bf Remark 2.} Cocylcles $\tau_2$ and $d \tau_2$ trivialize on the 
subalgebra
of divergence zero vector fields, because the divergence of a vector field is equal 
to the trace of the jacobian: div$(\v) = \Tr (\v^J)$.

\

{\bf Remark 3.} In coordinates, cocycles $\tau_1, \tau_2$ are written as follows 
(cf. [BB]):
$$\tau_1 \left( v(x) {\d \over \d x_i}, w(x) {\d \over \d x_j} \right) =
{\d v \over \d x_j} \sum\limits_{k=1}^N {\d^2 w \over \d x_i \d x_k} dx_k ,$$
$$\tau_2 \left( v(x) {\d \over \d x_i}, w(x) {\d \over \d x_j} \right) =
{\d v \over \d x_i} \sum\limits_{k=1}^N {\d^2 w \over \d x_j \d x_k} dx_k .$$

\

{\bf Remark 4.} When $N=1$, cocycles $\tau_1$ and $\tau_2$ coincide (the forms
$\Tr(AB)$ and $\Tr(A) \Tr(B)$ are equal in $\gl_1$), and become the Virasoro 
cocycle on the Witt algebra $W_1$. Thus for $N>1$, these cocycles are two 
distinct abelian 
generalizations of the Virasoro cocycle.
   
\

\

{\bf 2. Realization of the Virasoro-Bott cocycle via the Heisenberg group.}

\

 The Virasoro-Bott cocycle on the group of diffeomorphisms of $S^1$ has a rather 
mysterious form [Bo]:
$$ (F, \al)(G, \be) = \left(FG, \al + \be + B(F,G)\right), $$
where $F, G \in \Diff(S^1)$, $\al, \be \in \R$ and
$$B(F,G) = \int\limits_0^1 \ln|F^\prime(G(t))| \; d \ln|G^\prime (t)|.\eqno{(2.1)}$$
The cocycle $B(F,G)$ seems to have more analytic flavour rather than algebraic, 
and looks quite different from other cocycles on the infinite dimensional groups,
e.g., affine groups. We are now going to elucidate the situation with the unusual
form of this cocycle. We will link the Virasoro-Bott cocycle with the Heisenberg
cocycle, which is defined in more algebraic terms. 
The infinite dimensional Heisenberg group is the affine group 
which corresponds to $\GL_1$, thus the Heisenberg
cocycle may be viewed as the simplest case of the affine group cocycle.

The link between $\Diff(S^1)$ and the loop group of $\GL_1$ has been seen in fact
in (1.2) in the previous section, where the jacobian $J$ maps the former 
group into 
the latter. We will recover the Virasoro-Bott cocycle by evaluating the Heisenberg
cocycle on the jacobians of the diffeomorphisms.

 Let us recall the construction of the infinite dimensional Heisenberg group. We will
in fact give two versions of the Heisenberg group, one algebraic, and another analytic.

In our first construction, we will directly exponentiate the infinite dimensional 
Heisenberg Lie algebra. This Lie algebra has the basis 
$\left\{ H_j, c | j\in\Z \right\}$ and the Lie bracket
$$ [H_n, H_m] = n \delta_{n,-m} c ,$$
and $c$ is central. 

By the direct exponentiation we will get the Heisenberg group:
$$\left\{ \exp(\al c) \exp(\sum_{j\leq 0} a_j H_j) 
\exp(\sum_{j>0} a_j H_j) \right\}, 
\eqno{(2.2)}$$
where $a_k \in\R$ with only finitely many non-zero $a_k$. 
Using the formula
$$\exp(A) \exp(B) = \exp([A,B]) \exp(B) \exp(A),$$
which holds when $[A,B]$ commutes with $A$ and $B$, we get the multiplication
law in the Heisenberg group:
$$ \exp(\al c) \exp(\sum_{j\leq 0} a_j H_j) \exp(\sum_{j>0} 
a_j H_j)
 \times \exp(\be c) \exp(\sum_{j \leq 0} b_j H_j) \exp(\sum_{j>0} 
b_j H_j) $$
$$=
\exp((\al + \be + \sum_{j>0} j a_j b_{-j})c)
\exp(\sum_{j\leq 0} (\al_j+\be_j) H_j) \exp(\sum_{j>0} (\al_j+\be_j) H_j) 
.$$

This version is an algebraic version of the Heisenberg group since only finite
products are allowed in (2.2). The group of diffeomorphisms on the other hand is
of analytic nature, so we would need to consider a certain completion of the
algebraic version of the Heisenberg group. The construction of the analytic 
version of the Heisenberg group will follow the general pattern of the affine
extensions of the loop groups [PS].
 
Consider the group $\Map(S^1, \R^*)$ of $\Ci$ loops in $\R^* = \GL_1(\R)$,
and its connected component of identity $\Map_0(S^1, \R^*) = \Map(S^1, \R_+^*)$.
The group structure in a loop group is given by pointwise multiplication.
In order to define a central extension of $\Map(S^1, \R_+^*)$, we need to consider
contractions of elements $f \in \Map(S^1, \R_+^*)$ to identity. 
The general scheme of affine extensions simplifies in our case since 
there exists a canonical contraction for each $f$. This is due to the fact that
the exponential map from $\R$ (the Lie algebra of $\R_+^*$) to $\R_+^*$ is
bijective. 

Indeed let $x = \ln (f)$, $x \in \Map(S^1,\R)$. Then 
$\tilde f(\tau) = \exp(\tau x), \tau\in [0,1]$ is a homotopy between identity and $f$.
This canonical contraction respects multiplication: $\widetilde {fg} = 
\tf \tg$
because $\R^*$ is commutative, and the usual properties of the logarithm and 
the exponential functions hold. In the general scheme of affine extensions one has to
consider all possible contractions of elements of a loop group to identity. In our
case this is not necessary because of the existence of a canonical contraction.
The affine cocycle is given by the formula [PS]:
$$C(f,g) = \int_{S^1} \int_{[0,1]} \tf^{-1} d\tf \wedge d\tg \, \tg^{-1} . $$
Here $\tf, \tg$ are viewed as functions of two variables: $t$ -- the parameter along
$S^1$ and $\tau$ -- the parameter of the homotopy. Let $x(t) = \ln (f(t))$, 
$y(t) = \ln (g(t))$. Then we can calculate $C(f,g)$ explicitly in the 
following way:
$$C(f,g) = \int_{S^1} \int_{[0,1]} {1 \over \tf \tg} 
\left( {\d \tf \over \d \tau} {\d \tg \over \d t} - 
{\d \tf \over \d t} {\d \tg \over \d \tau} \right) d\tau \wedge dt $$
$$ = \int_{S^1} \int_{[0,1]} \tau 
\left( x(t){dy \over dt}- {dx \over dt} y(t)
\right) d\tau \wedge dt.$$
Integrating in $\tau$, and then integrating the second term by parts with respect to $t$
we get
$$ C(f,g) = \int_{S^1} x(t) d y(t) .$$
We see that the formula for the Heisenberg cocycle becomes very similar to the
formula for the Virasoro-Bott cocycle:
$$ C(f,g) = \int_{S^1} \ln(f) d \ln(g) . \eqno{(2.3)}$$
Now $x(t)$ and $y(t)$ are periodic $\Ci$ functions. These can be represented
by Fourier series:
$$x(t) = {a_0 \over 2} + \sum_{j>0} a_j \cos(2\pi j t) + \sum_{j>0} a_{-j} \sin(2\pi j t),$$
$$y(t) = {b_0 \over 2} + \sum_{j>0} b_j \cos(2\pi j t) + \sum_{j>0} b_{-j} \sin(2\pi j t).$$
The functions $x(t)$, $y(t)$ are of class $\Ci$ if and only if the Fourier coefficients
satisfy the condition
$$\lim_{j\to \pm\infty} j^k a_j = \lim_{j\to \pm\infty} j^k b_j = 0 \hbox{\rm \ for all \ }
k \in\N . \eqno{(2.4)}$$  
Now by direct integration of (2.3) we can get the formula for the cocycle in terms of the 
Fourier coefficients:
$$C(f,g) = \pi \sum_{j\in\Z} j a_j b_{-j}  .$$
The analytic version of the Heisenberg group is the set of pairs $(f,\alpha)$, where 
\break
$f\in \Map(S^1,\R^*_+)$, $\al\in\R$, and the group law is
$$ (f,\alpha)(g,\beta) = (fg, \al + \be + C(f,g)) .$$

It is easy to see that this cocycle is equivalent to the cocycle from the first version of 
the Heisenberg group.  Indeed, the isomorphism between the latter analytic version and the 
completion of the algebraic version can be given as follows:
$$ \varphi \left( \exp\left({a_0 \over 2} + \sum_{j>0} a_j \cos(2\pi j t) + \sum_{j>0} a_{-j} 
\sin(2\pi j t)\right), \al \right) = $$
$$ \exp(({\al\over 2\pi} + {1\over 2} \sum_{j>0} j a_j a_{-j})c)
\exp(\sum_{j\leq 0} a_j H_j) \exp(\sum_{j>0} a_j H_j) .$$
The condition on coefficients $\{ a_j \}$ for the completion of the first version of the 
Heisenberg group should be taken the same as the $\Ci$ condition (2.4).

\

{\bf Remark 5.} We can lift the cocycle (2.3) to the group $\Map(S^1, \R^*) =$ 
$\Map (S^1, \R_+^*) \times \Z_2$ by taking absolute values of the loops:
$$ C(f,g) = \int_{S^1} \ln|f| d \ln|g| . \eqno{(2.5)}$$

\

 Now we can give interpretation of the Virasoro-Bott cocycle by means of the cocycle of
the Heisenberg group.
  
 {\bf Theorem 2.1.} 
Let $J$ be the jacobian map $J: \Diff(S^1) \rightarrow \Map(S^1, \R^*)$.
The Virasoro-Bott cocycle on $\Diff(S^1)$ can be obtained by evaluating the 
Heisenberg cocycle (2.5) on $F^J (G)$ and $G^J$:
$$B(F,G) = C(F^J(G),G^J) .$$

Proof. This follows immediately from comparing the formulas (2.1) and (2.5).

\

{\bf Remark 6.} Using the jacobian map we also get a similar relation between the
Virasoro and Heisenberg cocycles at the level of Lie algebras. 
The Heisenberg cocycle on the abelian Lie algebra of functions on a circle is
$$ h( u_1 (x), u_2(x)) = \int_{S^1} u_1(x) d u_2 (x) .$$
The Virasoro cocycle
on $W_1$ is obtained by evaluating the Heisenberg cocycle on the jacobians 
(i.e. first derivatives) of vector fields:
$$ \left(v(x){\d\over \d x}, w(x){\d\over\d x}\right) = h(\v^J, \w^J) 
= \int_{S^1} v^\prime(x) d w^\prime(x) .$$

{\bf Remark 7.} Integration along $S^1$ in the formulas (2.3) and (2.5) realizes
the isomorphism between $\Omega^1 / d\Omega^0$ and $\R$. We thus can also write
the Heisenberg cocycle as a $\Omega^1 / d\Omega^0$ valued cocycle
$$ C(f,g) = \ln|f| d \ln|g| . \eqno{(2.6)}$$

\

{\bf Proposition 2.2.} Let $X$ be a $\Ci$ real manifold. Then the formula (2.6) gives 
an $\Omega^1(X) / d\Omega^0(X)$ valued central 2-cocycle on the group $\Map(X,\R^*)$. 
 
\

{\bf Proof.} We first verify the claim for the subgroup $\Map(X,\R^*_+)$. 
Applying the logarithm function, we get that the multiplicative group 
$\Map(X,\R^*_+)$ is isomorphic to the additive group of $\Ci$ functions on $X$.
Under this isomorphism (2.5) transforms into 
$$C_1 (u,v) = u(x) d v(x) , \eqno{(2.7)}$$
where $u = \ln(f), v = \ln(g)$. To see that (2.7) satisfies the cocycle condition, we 
need
to check that
$$ C_1 (u,v) + C_1(u+v,w) = C_1(v,w) + C_1(u, v+w) ,$$
which follows immediately from the linearity of (2.7). 

Finally, the group $\Map(X,\R^*)$ is a direct product $\Map(X,\R^*_+) \times \Z_2^n$,
where $n$ is the number of connected components of $X$,
thus a central cocycle on the subgroup $\Map(X,\R^*_+)$ lifts to a cocycle on
$\Map(X,\R^*)$.

\

Even though the above theorems are just  very simple observations, 
they allow us to construct  interesting generalizations of 
the Virasoro-Bott cocycle that will be explored in the next two sections.

\

\

{\bf 3. Abelian extension of the group of diffeomorphisms of a torus.}

\

The Virasoro cocycle can not be generalized as a central cocycle to the case of
a Lie algebra of vector fields in dimension greater than one. However in the case of a torus,
it is known that there are abelian generalizations of the Virasoro cocycle. We thus may
expect that we can also construct abelian extensions of the group of diffeomorphisms of a torus.

In this section we will construct an abelian extension that corresponds to Lie algebra cocycle
$d\tau_1$. The space of the extension will be $\Omega^2(\TN)$. This cocycle will be non-vanishing
on the subgroup of volume-preserving diffeomorphisms.

As we have seen in Section 1, the two abelian cocycles on $W_N$, $\tau_1$ and $\tau_2$, 
correspond to the decomposition of $\glN$ into $\slN \oplus \R I$. The cocycle that
will be constructed in this section corresponds to the $\slN$ component. In Section 4 we will
construct a cocycle that corresponds to the component of scalar matrices.  

Motivated by Theorem 2.1, our idea is to use the crossed homomorphism (1.2) and 
pull back a central 2-cocycle on $\Map(\TN, \GLN)$ to an abelian
2-cocycle on $\Diff(\TN)$ .

We will first describe our construction in a very general set-up.

Let a group $D$ acts on a group $M$ on the right by automorphisms:
$$ m \mapsto m(d), \hbox{\rm \ for \ } m\in M, d\in D .$$

{\bf Definition.} A crossed homomorphism is a map 
$$ j: D \rightarrow M,$$
such that
$$ j(d_1 d_2) = j(d_1) (d_2) \cdot j(d_2)  \eqno{(3.1)}$$
for all $d_1, d_2 \in D$.

Suppose we indeed have a crossed homomorphism $j: D \rightarrow M$.

Let $K$ be an abelian (additive) group and let 
$$ c: M \times M \rightarrow K$$
be a central 2-cocycle on $M$. Suppose that $D$ also acts on $K$
on the right by automorphisms.

{\bf Lemma 3.1.} If the $K$ valued 2-cocycle $c$ on $M$ is $D$ invariant, i.e., 
$$ c(m_1, m_2) (d) = c(m_1(d), m_2(d)), \hbox{\rm \ for all \ } m_1, m_2 \in M, d \in D,$$
then 
$$ b(d_1, d_2) = c(j(d_1)(d_2), j(d_2)) \eqno{(3.2)}$$ 
defines an abelian $K$ valued 2-cocycle on $D$. 

{Proof.} Using the fact that $c$ satisfies the central cocycle condition
$$c(m_1, m_2) + c(m_1 m_2, m_3) = c(m_1, m_2 m_3) + c(m_2, m_3) , \eqno{(3.3)}$$
we need to show that $b$ is an abelian cocycle:
$$b(d_1, d_2)(d_3) + b(d_1 d_2, d_3) = b(d_1, d_2 d_3) + b(d_2, d_3) . \eqno{(3.4)}$$
From the definition (3.2), we get that the left hand side is
$$c(j(d_1)(d_2), j(d_2)) (d_3) + c(j(d_1 d_2)(d_3), j(d_3)),$$
which using the invariance of cocycle $c$ and (3.1), becomes
$$c(j(d_1) (d_2 d_3), j(d_2)(d_3)) + c(j(d_1) (d_2 d_3) j(d_2)(d_3), j(d_3)) .$$   
By (3.3) this equals
$$c(j(d_1) (d_2 d_3), j(d_2)(d_3)j(d_3)) + c(j(d_2)(d_3), j(d_3)) $$
$$ = c(j(d_1) (d_2 d_3), j(d_2 d_3)) + c(j(d_2)(d_3), j(d_3)) $$
$$ = b(d_1, d_2 d_3) + b(d_2, d_3) ,$$
which is the left hand side of (3.4). The Lemma is proved.

We will apply this Lemma to the groups $D = \Diff(\TN)$ and $M = \Map(\TN, \GLN)$.
The extension group $K$ will be $\Omega^2 (\TN)$ and the crossed homomorphism between 
$\Diff(\TN)$ and $\Map(\TN, \GLN)$ is the jacobian map (1.2). The action of
$\Diff(\TN)$ on $\Map(\TN, \GLN)$ and $\Omega^2$ is given by (1.4) and 
(1.1).

The central 2-cocycle on $\Map(\TN, \GLN)$ that we will use is (cf. [PS], [LNMS]):
$$ C(f(\x),g(\x)) = \Tr \left( f^{-1} df \wedge dg \, g^{-1} \right). \eqno{(3.5)} $$

Let us verify the cocycle condition
$$C(f,g) + C(fg,h) = C(f, gh) + C(g,h). \eqno{(3.6)}$$
The left hand side here is
$$\Tr \left( f^{-1} df \wedge dg \, g^{-1} \right) + \Tr \left( (fg)^{-1} d(fg) \wedge dh \, h^{-1} \right) $$
$$ = \Tr \left( f^{-1} df \wedge dg \, g^{-1} \right) + 
\Tr \left( g^{-1} f^{-1} d(f)g  \wedge dh \, h^{-1} \right) +
\Tr \left( g^{-1} dg  \wedge dh \, h^{-1} \right) .$$
Similarly, the left hand side of (3.6) becomes
$$\Tr \left( f^{-1} df \wedge d(gh) (gh)^{-1} \right) + \Tr \left( g^{-1} dg \wedge dh \, h^{-1} \right) $$
$$ = \Tr \left( f^{-1} df \wedge dg \, g^{-1} \right) + 
\Tr \left( f^{-1} d(f)  \wedge  g dh \, h^{-1} g^{-1} \right) +
\Tr \left( g^{-1} dg  \wedge dh \, h^{-1} \right) .$$
Both sides are equal because
$$\Tr \left( g^{-1} f^{-1} d(f)g  \wedge dh \, h^{-1} \right) =
\Tr \left( f^{-1} d(f)  \wedge  g \, dh \, h^{-1} g^{-1} \right) $$
by the commutativity of the trace: $\Tr(AB) = \Tr (BA)$.

Finally we note that the cocycle (3.5) on $\Map(\TN, \GLN)$ is $\Diff(\TN)$ invariant,
since the action of diffeomorphisms on differential forms (1.1) agrees with 
the action on functions (1.4):
$$\Tr \left( f^{-1} df \wedge dg g^{-1} \right) (H) =
\Tr \left( f^{-1}(H) df(H) \wedge dg(H) g^{-1}(H) \right) .$$

We thus obtain the following result:

{\bf Theorem 3.2.} The set $\Diff(\TN) \times \Omega^2(\TN)$ becomes an extension
of the group of diffeomorphisms with the abelian group $\Omega^2(\TN)$ with multiplication
defined as follows:
$$(F,\al) (G,\be) = \left( FG, \al(G) + \be + \Tr \left( f^{-1} df \wedge dg g^{-1} \right)
\right),$$
where $f = F^J (G), g = G^J$.

\

{\bf Remark 8.} It is easy to compute that the corresponding Lie algebra extension
of $W_N$ is given by cocycle $d\tau_1$.

\

{\bf Problem.} Construct an abelian extension of $\Diff(\TN)$ with abelian subgroup
$\Omega^1 / d\Omega^0$ which corresponds to the Lie algebra 2-cocycle $\tau_1$.

\

\

{\bf 4. Abelian extension of the group of diffeomorphisms of a manifold with a volume form.}

\

In this section we introduce an abelian extension of the group of diffeomorphisms of
a manifold with a volume form that generalizes the Virasoro-Bott central extension.

Let $X$ be a real $\Ci$ manifold with the volume form $\omega$. We note that for
an arbitrary diffeomorphism $F\in\Diff(X)$, the quotient $\omega(F)\over \omega$ is
a well-defined non-vanishing function. Moreover, the map
$$ \delta: \Diff(X) \rightarrow \Map(X, \R^*),$$
given by
$$ \delta (F) = {\omega(F) \over \omega}$$
is a crossed homomorphism. Indeed,
$$\delta(FG) = {\omega(FG) \over \omega} = {\omega(FG) \over \omega(G)} {\omega(G) \over \omega}
= \delta(F) (G) \delta(G) .$$ 

As an immediate consequence of this observation, Lemma 3.1 and Proposition 2.2, 
we get the following:

{\bf Theorem 4.1.} Let $X$ be a $\Ci$ real manifold with a volume form $\omega$ and
let $C(f,g)$ be the Heisenberg cocycle (2.6) on $\Map(X,\R^*)$.
Then
$$ B(F,G) = C \left({\omega(FG) \over \omega(G)}, {\omega(G) \over \omega} \right) $$
$$ = \ln \left|{\omega(FG) \over \omega(G)}\right| d\ln \left| {\omega(G) \over \omega} \right|
\eqno{(4.1)}$$
is an abelian $\Omega^1(X) / d\Omega^0(X)$ valued 2-cocycle on the group of diffeomorphisms
$\Diff(X)$.

\

{\bf Corollary 4.2.} The 2-cocycle on the Lie algebra of vector fields on
a manifold with a volume form $\omega$ that corresponds to the group cocycle (4.1) is
$$ c(\v,\w) = \diver(\v) \, d \, \diver(\w),$$
where the divergence of a vector field $\v$ is given by 
$\diver(\v) = {\v \cdot \omega \over \omega}$.

\

{\bf Example 4.3.} If $X$ is a torus $\TN$ with the volume form $\omega = dx_1 \wedge \ldots
\wedge dx_N$, then ${\omega(F) \over \omega} = det(F^J)$, and the cocycle (3.1) can be written
as
$$B(F,G) = \ln |det(F^J(G))| d \ln |det(G^J)| .$$
One can easily calculate that the corresponding Lie algebra cocycle on $W_N$ is $\tau_2$.

\

{\bf Acknowledgments.} I am grateful to l'Institut des Math\'ematiques
de Jussieu (Paris) and Mathematical Sciences Research Institute
(Berkeley) for their hospitality. I also thank Dmitri Orlov for stimulating
discussions.

\

{\bf Bibliography.}

\

\noindent
[BB] {Berman, S. and Billig, Y.}:
{Irreducible Representations for Toroidal Lie algebras,}
{\it J. Algebra}, {\bf 221} (1999), 188--231.

\noindent
[Bi] {Billig, Y.}: 
{Energy-momentum tensor for the toroidal Lie algebras},
preprint, 
\hfill\break
math.RT/0201313.

\noindent
[Bo] {Bott, R.}: {On the characteristic classes of groups of diffeomorphisms},
{\it Enseign. Math. (2)}, {\bf 23} (1977), 209--220.

\noindent
[D]  {Dzhumadil'daev, A.}: 
{Virasoro type Lie algebras and deformations,} 
{\it Z. Phys. C}, {\bf 72} (1996), 509--517.
 
\noindent
[EM] {Eswara Rao, S. and Moody, R.V.}:
{Vertex Representations for $N$-Toroidal Lie Algebras and a 
Generalization of the Virasoro Algebra,}
{\it Comm. Math. Phys.}, {\bf 159}(1994), 239--264.

\noindent
[FK] {Frenkel, I. and Khesin, B.}: 
{Four-dimensional realization of two-{}dimensional current groups}, 
{\it Comm. Math. Phys.}, {\bf 178} (1996), 541--562.

\noindent
[K] {Kurosh, A.G.}: The theory of groups, Chelsea, N.Y., 1956. 

\noindent
[L] {Larsson, T.A.}:
{Lowest-Energy Representations of Non-Centrally Extended Diffeomorphism Algebras,}
{\it Comm. Math. Phys.}, {\bf 201} (1999), 461--470.

\noindent
[LMNS] {Losev, A., Moore, G., Nekrasov, N. and Shatashvili, S.}: 
{Central extensions of gauge groups revisited}, 
{\it Selecta Math. (N.S.)}, {\bf 4}(1) (1998), 117--123.

\noindent
[OR] {Ovsienko, V. and Roger, C.}: 
{Generalizations of Virasoro group and Virasoro
algebra through extensions by modules of tensor densities on $S^1$},
{\it Indag. Math. N.S.}, {\bf 9} (2), (1998), 277--288.

\noindent
[PS] {Pressley, A. and  Segal, G.}: Loop groups, Oxford University Press, Oxford, 1986.

\end